# Distributed Model Predictive Control for Heterogeneous Vehicle Platoons under Unidirectional Topologies

Yang Zheng, Shengbo Eben Li, Keqiang Li, Francesco Borrelli, and J. Karl Hedrick

*Abstract*— This paper presents a distributed model predictive control (DMPC) algorithm for heterogeneous vehicle platoons with unidirectional topologies and *a priori* unknown desired set point. The vehicles (or nodes) in a platoon are dynamically decoupled but constrained by spatial geometry. Each node is assigned a local open-loop optimal control problem only relying on the information of neighboring nodes, in which the cost function is designed by penalizing on the errors between predicted and assumed trajectories. Together with this penalization, an equality based terminal constraint is proposed to ensure stability, which enforces the terminal states of each node in the predictive horizon equal to the average of its neighboring states. By using the sum of local cost functions as a Lyapunov candidate, it is proved that asymptotic stability of such a DMPC can be achieved through an explicit sufficient condition on the weights of the cost functions. Simulations with passenger cars demonstrate the effectiveness of proposed DMPC.

*Index Terms*— Autonomous vehicle, heterogeneous platoon, graph theory, distributed control, model predictive control

## I. Introduction

THE platooning of autonomous vehicles has received considerable attention in recent years [1]-[7]. Most of this attention is due to its potential to significantly benefit road transportation, including improving traffic efficiency, enhancing road safety and reducing fuel consumption, *etc.* [1][2]. The main objective of platoon control is to ensure all the vehicles in a group move at the same speed while maintaining a pre-specified distance between any consecutive followers [5]-[7].

The earliest practices on platoon control could date back to the PATH program in the eighties of the last century, in which many well-known topics were introduced in terms of sensors and actuators, control architecture, decentralized control and string stability [5]. Since then, many other issues on platoon control have been discussed, such as the selection of spacing policies [6][7], influence of communication topology [8]-[10], and impact of dynamic heterogeneity [11][12]. In recent years, some advanced platoon control laws have been proposed under the framework of multi-agent consensus control, see *e.g.* [13]-[17]. Most of them employ linear dynamics and linear controllers for the convenience of theoretical completeness, and do not account for input constraints and model nonlinearities. A few notable exceptions are in [14] and [17], where the communication topologies are assumed to be limited in range. However, the input constraints and model nonlinearities do exist in a more accurate problem formulation due to actuator saturation and some salient nonlinearities involved in the powertrain system, *e.g.*, engine, driveline and aerodynamic drag [18][19]. Besides, with the rapid deployment of vehicle-to-vehicle (V2V) communication, such as DSRC and VANETs [20], various types of communication topologies are emerging, *e.g.,* the two-predecessor following type and the multiple-predecessor following type [21][22]. New challenges for platoon control arise naturally considering the variety of topologies, especially when taking into account a large variety of topologies in a systematic and integrated way.

This paper proposes an innovative solution for platoon control considering both nonlinear dynamics and topological variety based on the model predictive control (MPC) framework. Traditionally, MPC is used for a single-agent system, where the control input is obtained by numerically optimizing a finite horizon optimal control problem where both nonlinearity and constraints can be explicitly handled [23]. This technique has been embraced by many industrial applications, for instance, thermal energy control [24], collision avoidance [18], and vehicle stability [25], and energy management [26], *etc*. Most of these MPCs are implemented in a centralized way, where all the control inputs are computed by assuming all the states are known, *e.g.*, [18][24]-[26]. When considering an actual platoon system involving multiple vehicles, the centralized implementation is not suitable because of the limitation to gather the information of all vehicles and the challenge to compute a large-scale optimization problem. In this paper, we present a synthesis method of distributed model predictive control (DMPC) for a heterogeneous platoon, where each vehicle is assigned a local optimal control problem only relying on its neighboring vehicles' information.

Recently, several DMPC schemes have been proposed for dynamically coupled or decoupled multi-agent systems [27]-[30]. The asymptotic stability was usually established by

This research is supported by NSF China with grant 51575293 and CSC funding with grant 201406215042.

Y. Zheng was with the State Key Lab of Automotive Safety and Energy, Tsinghua University, Beijing 100084, China. He is now with the Department of Engineering Science, University of Oxford, Parks Road, Oxford OX1 3PJ, U.K. (e-mail: yang.zheng@eng.ox.ac.uk)

S Eben Li is with the State Key Lab of Automotive Safety and Energy, Tsinghua University, Beijing 100084, China. He is currently working on University of California, Berkeley as a visiting scholar. (Corresponding author, e-mail: lisb04@gmail.com)

K. Li is with the State Key Lab of Automotive Safety and Energy, Tsinghua University, Beijing 100084, China. (e-mail: likq@tsinghua.edu.cn)

F. Borrelli, and J. K. Hedrick are with the Department of Mechanical Engineering, University of California, Berkeley, CA 94720 USA. (e-mail: fborrelli@berkeley.edu; karlhed@gmail.com)

employing the consistency constraints, *e.g.*, the mismatch between newly calculated optimal trajectories and the previously calculated ones must be bounded [27][28]. A recent comprehensive review on DMPC can be found in [31]. However, the majority of existing DMPC algorithms only focus on the stabilization of the system with a common set point, assuming all agents *a priori* know the desired equilibrium information. For a vehicle platoon, such a common set point corresponds to the leader's state. However, it is not practical to assume all the followers can communicate with the leader, which means not all of the followers know the desired set point in a platoon. The purpose of this paper is to address the control issue of vehicle platoons with *a priori* unknown desired set point under distributed MPC framework. Most existing MPC works in this field rely on the problem formulation of adaptive cruise control (ACC), *e.g.*, [32][33], which only involve two vehicles in the problem formulation. There exist some extensions to the cooperative ACC case which involve multiple vehicles, *e.g.*, [14][34]. Such treatments in [14] and [34], however, also directly take two consecutive vehicles into the problem formulation, which are only applicable to limited types of communication topologies, *i.e.*, the predecessor-following type and the predecessor-leader following type.

This paper presents a distributed model predictive control (DMPC) algorithm for heterogeneous platoons with unidirectional topologies and *a priori* unknown desired set point. The contribution of this paper is in two aspects: 1) the proposed DMPC algorithm does not need all nodes to *a priori* know the desired set point, which is a significant improvement compared to many previous studies, *e.g.*. [27]-[30]; 2) our finding not only explicitly highlights the importance of communication topology in stabilizing the entire platoon system, but also extends the results in [14][34] to suit any arbitrary unidirectional topology. Specifically, a platoon is viewed as a group of vehicles, which are dynamically decoupled but interact with each other by spatial geometry and communication topology. In a platoon, only the followers, which directly communicate with the leader, know the desired set point. Under the proposed DMPC, each follower is assigned a local open-loop optimal control problem only relying on the information of neighboring vehicles, in which the errors between predicted trajectories and assumed ones are penalized. A neighboring average based terminal constraint is proposed, by which the terminal state of each node is enforced to be equal to the state average of its neighboring nodes. We use the sum of the local cost functions as a Lyapunov candidate, and prove that asymptotical stability can be achieved through explicit parametric conditions on the weights of the cost functions under unidirectional topologies. The material in this paper was partially summarized in [2].

The rest of this paper is organized as follows. In section II, the dynamic model, control objective and model of communication topology in a platoon are presented. Section III introduces the formulation of local optimal control problems. The stability results are given in section IV, followed by the simulation results in section V. Section VI concludes this paper.

*Notation:* Throughout this paper, $\mathbb{R}, \mathbb{C}$ stand for the set of real numbers and complex numbers, respectively. We use $\mathbb{R}^{m \times n}$ to denote the set of $m \times n$ real matrices, and the set of symmetric matrices of order $n$ is denoted by $\mathbb{S}^n$. For any positive integer $N$, let $\mathcal{N} = \{1,2,\cdots,N\}$. Given a symmetric matrix $M \in \mathbb{S}^n$, $M \geq 0$ ($M > 0$) means that the matrix is positive semi-definite (positive definite). The relation $M_1 \geq M_2$ for symmetric matrices means that $M_1 - M_2 \geq 0$. The identity matrix of dimension $n$ is denoted by $I_n$. $\text{diag}(a_1,\cdots,a_N)$ is a diagonal matrix with main diagonal entries $a_j, j \in \mathcal{N}$ and the off-diagonal entries are zero. Given a matrix $A \in \mathbb{R}^{n \times n}$, its spectrum radius is denoted by $\rho(A)$. Given a vector $x$ and a positive semi-definite matrix $Q \geq 0$, we use $\|x\|_Q = (x^T Q x)^{1/2}$ to denote the weighted Euclidean norm. The Kronecker product is denoted by $\otimes$, which facilitates the manipulation of matrices by the following properties 1)$(A \otimes B)(C \otimes D) = AC \otimes BD$, 2)$(A \otimes B)^T = A^T \otimes B^T$.

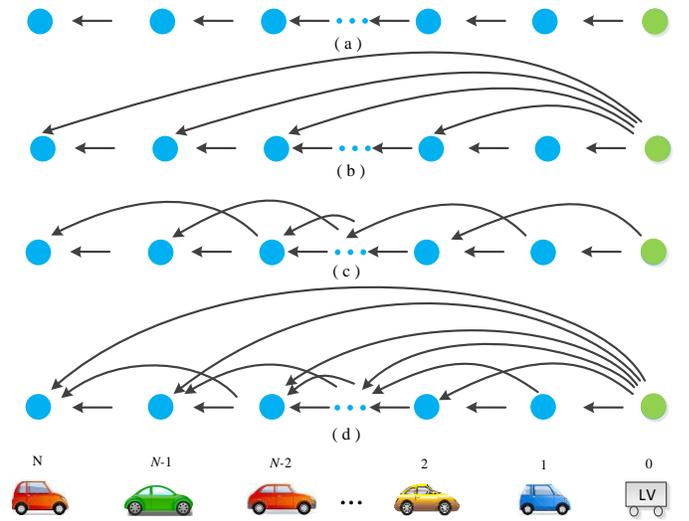

Fig. 1. Examples of unidirectional topology: (a) predecessor-following (PF), (b) predecessor-leader following (PLF), (c) two-predecessor following (TPF), (d) two-predecessor-leader following (TPLF).

## II. PLATOON MODELING AND CONTROL OBJECTIVE

As shown in Fig. 1, this paper considers a heterogeneous platoon with a broad selection of communication topologies running on a flat road with $N + 1$ vehicles (or nodes), which includes a leading vehicle (LV, indexed by 0) and $N$ following vehicles (FVs, indexed from 1 to $N$). The communication among nodes is assumed to be unidirectional from the preceding vehicles to downstream ones, which are commonly used in the field of vehicle platoon [22], such as predecessor-following (PF), predecessor-leader following (PLF), two-predecessor following (TPF), and two-predecessor-leader following (TPLF) (see Fig. 1 for examples).

The platoon is dynamically decoupled, but constrained by the spatial formation. Each node has nonlinear dynamics with input constraints, but its desired set point with respect to the leader might be unknown. Only the nodes that directly communicate with the leader know the desired set point. The control objective of DMPC is to achieve a global coordination in terms of movement and geometry even though the exchanged information is local and limited to the neighborhood of each node.

### A. Nonlinear platoon model for control

This paper only considers the vehicle longitudinal dynamics, which are composed of engine, drive line, brake system, aerodynamic drag, tire friction, rolling resistance, gravitational force, *etc*. To strike a balance between accuracy and conciseness, it is assumed that: (1) the vehicle body is rigid and left-right symmetric; (2) the platoon is on flat and dry-asphalt road, and the tyre slip in the longitudinal direction is neglected; (3) the powertrain dynamics are lumped to be a first-order inertial transfer function; (4) the driving and braking torques are integrated into one control input. Then, the discrete-time model of any following vehicle $i$ is

$$\begin{cases} s_i(t+1) = s_i(t) + v_i(t)\Delta t \\ v_i(t+1) = v_i(t) + \dfrac{\Delta t}{m_{\text{veh},i}}\left(\dfrac{\eta_{\text{T},i}}{R_i}T_i(t) - F_{\text{veh},i}(v_i(t))\right) \\ T_i(t+1) = T_i(t) - \dfrac{1}{\tau_i}T_i(t)\Delta t + \dfrac{1}{\tau_i}u_i(t)\Delta t \end{cases}$$

$$F_{\text{veh},i}(v_i(t)) = C_{A,i}v_i^2(t) + m_{\text{veh},i}gf_i \qquad (1)$$

where $\Delta t$ is the discrete time interval; $s_i(t), v_i(t)$ denote the position and velocity of node $i$; $m_{\text{veh},i}$ is the vehicle mass; $C_{A,i}$ is the coefficient of aerodynamic drag; $g$ is the gravity constant; $f_i$ is the coefficient of rolling resistance; $T_i(t)$ is the integrated driving/braking torque; $\tau_i$ is the inertial lag of longitudinal dynamics, $R_i$ is the tire radius; $\eta_{\text{T},i}$ is the mechanical efficiency of the driveline, and $u_i(t) \in \mathbb{R}$ is the control input, representing the desired driving/braking torque. The control input is subject to the box constraint:

$$u_i \in \mathcal{U}_i = \{u_{\min,i} \leq u_i \leq u_{\max,i}\}, \qquad (2)$$

where $u_{\min,i}, u_{\max,i}$ are the bounds. For each node, the state is denoted as $x_i(t) = [s_i(t), v_i(t), T_i(t)]^T \in \mathbb{R}^{3\times 1}$, and the output is denoted as $y_i(t) = [s_i(t), v_i(t)]^T \in \mathbb{R}^{2\times 1}$. Further, (1) can be rewritten into a compact form

$$\begin{aligned} x_i(t+1) &= \phi_i(x_i(t)) + \psi_i \cdot u_i(t), \\ y_i(t) &= \gamma x_i(t) \end{aligned} \qquad (3)$$

where $\psi_i = \left[0, 0, \dfrac{1}{\tau_i}\Delta t\right]^T \in \mathbb{R}^{3\times 1}$, $\gamma = \begin{bmatrix} 1 & 0 & 0 \\ 0 & 1 & 0 \end{bmatrix} \in \mathbb{R}^{2\times 3}$, $\phi_i(x_i) \in \mathbb{R}^{3\times 1}$ is defined as

$$\phi_i = \begin{bmatrix} s_i(t) + v_i(t)\Delta t \\ v_i(t) + \dfrac{\Delta t}{m_{\text{veh},i}}\left(\dfrac{\eta_{\text{T},i}}{R_i}T_i(t) - F_{\text{veh},i}(v_i(t))\right) \\ T_i(t) - \dfrac{1}{\tau_i}T_i(t)\Delta t \end{bmatrix}.$$

Define $X(t) \in \mathbb{R}^{3N\times 1}$, $Y(t) \in \mathbb{R}^{2N\times 1}$, and $U(t) \in \mathbb{R}^{N\times 1}$ as the vectors of states, outputs and inputs of all nodes, *i.e.*,

$$X(t) = [x_1^T(t), x_2^T(t), \cdots, x_N^T(t)]^T,$$
$$Y(t) = [y_1^T(t), y_2^T(t), \cdots, y_N^T(t)]^T,$$
$$U(t) = [u_1(t), \cdots, u_N(t)]^T.$$

Then, the overall discrete-time dynamics of the platoon becomes:

$$\begin{aligned} X(t+1) &= \mathbf{\Phi}(X(t)) + \mathbf{\Psi} \cdot U(t), \\ Y(t+1) &= \mathbf{\Gamma} \cdot X(t+1), \end{aligned} \qquad (4)$$

where $\mathbf{\Phi} = [\phi_1(x_1)^T, \phi_2(x_2)^T, \cdots, \phi_N(x_N)^T]^T \in \mathbb{R}^{3N\times 1}$, $\mathbf{\Psi} = \text{diag}\{\psi_1, \cdots, \psi_N\} \in \mathbb{R}^{3N\times N}$, $\mathbf{\Gamma} = I_N \otimes \gamma \in \mathbb{R}^{2N\times 3N}$.

The model (1) for vehicle dynamics is inherently a third-order nonlinear system, which can encapsulate a wide range of vehicles. Note that linear models are also widely used in platoon control for the sake of theoretical completeness, *e.g.*, [6][13][15].

### B. Objective of platoon control

The objective of platoon control is to track the speed of the leader while maintaining a desired gap between any consecutive vehicles which is specified by a desired spacing policy, *i.e.*,

$$\begin{cases} \lim_{t\to\infty}\|v_i(t) - v_0(t)\| = 0 \\ \lim_{t\to\infty}\|s_{i-1}(t) - s_i(t) - d_{i-1,i}\| = 0 \end{cases}, i \in \mathcal{N}, \qquad (5)$$

where $d_{i-1,i}$ is the desired space between $i-1$ and $i$. The selection of $d_{i-1,i}$ determines the geometry formation of the platoon. Here, the constant spacing policy is used, *i.e.*,

$$d_{i-1,i} = d_0. \qquad (6)$$

### C. Model of communication topology

An accurate model of topological structure is critical to design a coupled cost function in DMPC. The communication topology in a platoon can be modeled by a directed graph $\mathbb{G} = \{\mathbb{V}, \mathbb{E}\}$, where $\mathbb{V} = \{0,1,2,\ldots,N\}$ is the set of nodes and $\mathbb{E} \subseteq \mathbb{V} \times \mathbb{V}$ is the set of edges in connection [21][36]. The properties of graph $\mathbb{G}$ is further reduced into the formation of three matrices, *i.e.*, adjacency matrix $\mathcal{A}$, Laplacian matrix $\mathcal{L}$ and pinning matrix $\mathcal{P}$.

The adjacency matrix is used to describe the directional communication among the followers, which is defined as $\mathcal{A} = [a_{ij}] \in \mathbb{R}^{N\times N}$ with each entry expressed as

$$\begin{cases} a_{ij} = 1, & if\ \{j,i\} \in \mathbb{E} \\ a_{ij} = 0, & if\ \{j,i\} \notin \mathbb{E} \end{cases}, i, j \in \mathcal{N}, \qquad (7)$$

where $\{j, i\} \in \mathbb{E}$ means there is a directional edge from node $j$ to node $i$, *i.e.*, node $i$ can receive the information of node $j$ (or simply $j \to i$). The Laplacian matrix $\mathcal{L} \in \mathbb{R}^{N\times N}$ is defined as

$$\mathcal{L} = \mathcal{D} - \mathcal{A}. \qquad (8)$$

where $\mathcal{D} \in \mathbb{R}^{N\times N}$ is called the in-degree matrix, defined as

$$\mathcal{D} = \text{diag}\{\deg_1, \deg_2, \cdots, \deg_N\}, \qquad (9)$$

where $\deg_i = \sum_{j=1}^{N} a_{ij}$ represents the in-degree of node $i$ in $\mathbb{G}$. The pinning matrix $\mathcal{P} \in \mathbb{R}^{N\times N}$ is used to model how each follower connects to the leader, defined as

$$\mathcal{P} = \text{diag}\{p_1, p_2, \ldots, p_N\}, \qquad (10)$$

where $p_i = 1$ if edge $\{0, i\} \in \mathbb{E}$; otherwise, $p_i = 0$. Node $i$ is said to be pinned to the leader if $p_i = 1$, and only the nodes pinned to the leader know the desired set point. We further define the leader accessible set of node $i$ as

$$\mathbb{P}_i = \begin{cases} \{0\}, & if\ p_i = 1 \\ \emptyset, & if\ p_i = 0 \end{cases}.$$

For the sake of completeness, several definitions are stated as

follows:

**1) Directed path.** A directed path from node $i_1$ to node $i_k$ is a sequence of edges $(i_1, i_2), (i_2, i_3), \cdots, (i_{k-1}, i_k)$ with $(i_{j-1}, i_j) \in \mathbb{E}, \forall j = \{2, \cdots, k\}$

**2) Spanning tree.** The graph $\mathbb{G}$ is said to contain a spanning tree if there is a root node such that there exists a directed path from this node to every other node.

**3) Neighbor set.** Node $j$ is said to be a neighbor of node $i$ if and only if $a_{ij} = 1, j \in \mathcal{N}$. The neighbor set of node $i$ is denoted by $\mathbb{N}_i = \{j | a_{ij} = 1, j \in \mathcal{N}\}$.

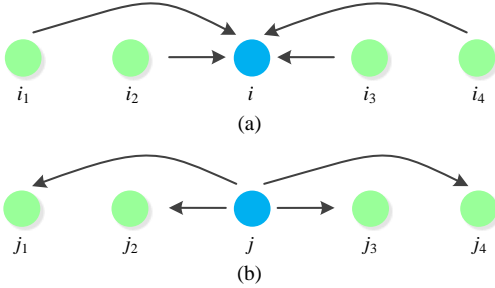

Fig. 2. Examples of sets $\mathbb{N}_i$ and $\mathbb{O}_i$: (a) $\mathbb{N}_i = \{i_1, i_2, i_3, i_4\}$, (b) $\mathbb{O}_j = \{j_1, j_2, j_3, j_4\}$.

The set $\mathbb{N}_i$ means that node $i$ can receive the information of any $j \in \mathbb{N}_i$. Similarly, we define a dual set $\mathbb{O}_i = \{j | a_{ji} = 1, j \in \mathcal{N}\}$, which means that node $i$ sends its information to any $j \in \mathbb{O}_i$. Note that for an undirected topology, we have $\mathbb{N}_i = \mathbb{O}_i$; but for any directed topology, this equality does not hold. Fig. 2 illustrates typical examples of sets $\mathbb{N}_i$ and $\mathbb{O}_i$.

Note that the set $\mathbb{I}_i = \mathbb{N}_i \cup \mathbb{P}_i$ describes all nodes which can send their information to node $i$. Hence, only the information of nodes in $\mathbb{I}_i$ can be used to construct the local optimal control problem for node $i$.

## III. DESIGN OF DISTRIBUTED MODEL PREDICTIVE CONTROL

This section introduces the formulation of DMPC for heterogeneous platoons. The position and velocity of the leader are denoted by $s_0(t)$ and $v_0(t)$ respectively. The leader is assumed to run at a constant speed, i.e., $s_0 = v_0 t$. The desired set point of state and input of node $i$ is

$$\begin{cases} x_{\text{des},i}(t) = [s_{\text{des},i}(t), v_{\text{des},i}(t), T_{\text{des},i}(t)]^T \\ u_{\text{des},i}(t) = T_{\text{des},i}(t), \end{cases} \quad (11)$$

where $s_{\text{des},i}(t) = s_0(t) - i \cdot d_0, v_{\text{des},i}(t) = v_0$ and $T_{\text{des},i}(t) = h_i(v_0)$, which is used to counterbalance the external drag, defined as

$$h_i(v_0) = \frac{R_i}{\eta_{T,i}} (C_{A,i} v_0^2 + m_{\text{veh},i} g f_i). \quad (12)$$

The corresponding equilibrium of output is $y_{\text{des},i}(t) = \gamma x_{\text{des},i}(t)$. Note that the constant speed assumption for the leader characterizes the desired equilibrium for a platoon, which is widely used for theatrical analysis in the literature [3][9][13]-[16]. Note also that many previous works on DMPC assume that all nodes *a priori* know the desired set point, *e.g.*, [14][27][28]. In this paper, it must be pointed out that the desired set point is not universally known for all followers in a platoon, and only the nodes pinned to the leader have access to the desired set information. The method proposed in this paper can guarantee the consensus of the desired set point among the followers when $\mathbb{G}$ contains a spanning tree.

### A. Local open-loop optimal control problem

For each node $i$, the formulation of its local optimal control problem only uses the information of the nodes in set $\mathbb{I}_i = \mathbb{N}_i \cup \mathbb{P}_i$. For the sake of narrative convenience, the nodes in $\mathbb{N}_i$ are numbered as $i_1, i_2, \cdots, i_m$. Define

$$y_{-i}(t) = [y_{i_1}^T(t), y_{i_2}^T(t), \cdots, y_{i_m}^T(t)]^T,$$
$$u_{-i}(t) = [u_{i_1}(t), u_{i_2}(t), \cdots, u_{i_m}(t)]^T$$

as the vectors of outputs and inputs of nodes in $\mathbb{N}_i$, respectively. The same length of predictive horizon $N_p$ is used in all local problems. Over the prediction horizon $[t, t + N_p]$, we define three types of trajectories:

$y_i^p(k|t)$: Predicted output trajectory,
$y_i^*(k|t)$: Optimal output trajectory,
$y_i^a(k|t)$: Assumed output trajectory,

where $k = 0, 1 \dots, N_p$. The notation $y_i^p(k|t)$ represents the output trajectory that parameterizes the local optimal control problem. The notation $y_i^*(k|t)$ represents the optimal solution after numerically solving the local problem. The notation $y_i^a(k|t)$ is the assumed trajectory transmitted to the nodes in set $\mathbb{O}_i$, which is actually the shifted last-step optimal trajectories of node $i$ (see the precise definition in Section III.B). Likewise, three types of control inputs are also defined,

$u_i^p(k|t)$: Predicted control input,
$u_i^*(k|t)$: Optimal control input,
$u_i^a(k|t)$: Assumed control input.

Now we define the local open-loop optimal control problem for each node $i$:

**Problem $\mathcal{F}_i$:** For $i \in \{1, 2, \dots, N\}$ at time $t$

$$\min_{u_i^p(0|t), \cdots, u_i^p(N_p-1|t)} J_i(y_i^p, u_i^p, y_i^a, y_{-i}^a)$$
$$= \sum_{k=0}^{N_p-1} l_i\left(y_i^p(k|t), u_i^p(k|t), y_i^a(k|t), y_{-i}^a(k|t)\right), \quad (13a)$$

subject to

$$x_i^p(k+1|t) = \phi_i\left(x_i^p(k|t)\right) + \psi_i \cdot u_i^p(k|t),$$
$$y_i^p(k|t) = \gamma \cdot x_i^p(k|t), \quad (13b)$$
$$x_i^p(0|t) = x_i(t),$$

$$u_i^p(k|t) \in \mathcal{U}_i, \quad (13c)$$

$$y_i^p(N_p|t) = \frac{1}{|\mathbb{I}_i|} \sum_{j \in \mathbb{I}_i} (y_j^a(N_p|t) + \tilde{d}_{i,j}), \quad (13d)$$

$$T_i^p(N_p|t) = h_i\left(v_i^p(N_p|t)\right), \quad (13e)$$

where $[u_i^p(0|t), \cdots, u_i^p(N_p - 1|t)]$ denotes the unknown variables to be optimized, $|\mathbb{I}_i|$ is the cardinality of set $\mathbb{I}_i$, and $\tilde{d}_{i,j} = [d_{i,j}, 0]^T$ denotes the desired distance vector between $i$ and $j$. The terminal constraint (13d) is to enforce that node $i$ has the same output as the average of assumed outputs in $\mathbb{I}_i$ at the end of predictive horizon. The terminal constraint (13e) is to enforce that node $i$ moves at constant speed without acceleration or deceleration at the end of predictive horizon. These two terminal constraints are critical to the stability of proposed DMPC.

The function $l_i$ in (13a) is the cost associated with node $i$, defined as

$$l_i\left(y_i^p(k|t), u_i^p(k|t), y_i^a(k|t), y_{-i}^a(k|t)\right)$$
$$= \left\|y_i^p(k|t) - y_{\text{des},i}(k|t)\right\|_{Q_i}$$
$$+ \left\|u_i^p(k|t) - h_i\left(v_i^p(k|t)\right)\right\|_{R_i} \quad (14)$$
$$+ \left\|y_i^p(k|t) - y_i^a(k|t)\right\|_{F_i}$$
$$+ \sum_{j \in \mathbb{N}_i} \left\|y_i^p(k|t) - y_j^a(k|t) - \tilde{d}_{i,j}\right\|_{G_i},$$

where $Q_i \in \mathbb{S}^2$, $R_i \in \mathbb{R}$, $F_i \in \mathbb{S}^2$ and $G_i \in \mathbb{S}^2$ are the weighting matrices. All weighting matrices are assumed to be symmetric and satisfy the following conditions:
a) $Q_i \geq 0$, which represents the strength to penalize the output error from the desired equilibrium. Note that $Q_i$ also contains the information whether node $i$ is pinned to the leader. If $p_i = 0$, node $i$ is unable to know its desired set point, and therefore $Q_i = 0$ is always enforced. If $p_i = 1$, then $Q_i > 0$ in its penalization functions.
b) $R_i \geq 0$, which represents the strength to penalize the input error diverged from equilibrium, meaning that the controller prefers to maintain constant speed.
c) $F_i \geq 0$, which means that node $i$ tries to maintain its assumed output. Note that this assumed output is actually the shifted last-step optimal trajectory of the same node, and this output is sent to the nodes in set $\mathbb{O}_i$.
d) $G_i \geq 0$, which means that node $i$ tries to maintain the output as close to the assumed trajectories of its neighbors (i.e., $j \in \mathbb{N}_i$) as possible.

**Remark 1**. The construction of (13d) is based on the local average of neighboring outputs, which is called neighboring average based terminal constraint. Thus, any node does not need to *a priori* know the desired set point, which must rely on pinning to the leader. This design is a significant improvement compared to many previous studies, which assumes that all nodes inherently pin to the leader if not explicitly mentioned, or only consider the stabilization of *a priori* known set point, *e.g.*, [27]-[30].

**Remark 2**. The formulation of problem $\mathcal{F}_i$ only needs the information from its neighbors, thus it is suitable for various communication topologies, including all of those shown in Fig. 1. However, stability might be not ensured by a normal DMPC law given by (13). A sufficient condition is needed to rigorously ensure asymptotic stability, which will be discussed and proved in section IV.

**Remark 3**. Note that problem $\mathcal{F}_i$ needs a precise vehicle model to predict the future output behavior. In addition to asymptotic stability (the focus of this paper), another challenging issue is the robustness to model uncertainty and noise, which is an active research topic [37]. Existing methods include robust optimization, worst-case and scenario-based approaches for constrained linear systems with disturbances [38]. We notice that a robust constraint proposed in [39] might be integrated into problem $\mathcal{F}_i$ to address the robustness issue. Besides, the coupling constraints for collision avoidance in multi-agent systems [18][30] are not considered in problem $\mathcal{F}_i$, which deserves further research.

### B. Algorithm of distributed model predictive control

The DMPC algorithm is shown as follows:
*(I) Initialization:*
At time $t = 0$, assume that all followers are moving at a constant speed, and initialize the assumed values for node $i$ as:

$$\begin{cases} u_i^a(k|0) = h_i(v_i(0)) \\ y_i^a(k|0) = y_i^p(k|0) \end{cases}, k = 0,1,\cdots,N_p - 1, \quad (15)$$

where $y_i^p$ is iteratively calculated by

$$x_i^p(k+1|0) = \phi_i\left(x_i^p(k|0)\right) + \psi_i \cdot u_i^a(k|0)$$
$$y_i^p(k|0) = \gamma\, x_i^p(k|0),\ x_i^p(0|0) = x_i(0).$$

*(II) Iteration of DMPC*:
At any time $t > 0$, for all node $i = 1, \cdots, N$,
1) Optimize Problem $\mathcal{F}_i$ according to its current state $x_i(t)$, its own assumed output $y_i^a(k|t)$, and assumed outputs from its neighbors $y_{-i}^a(k|t)$, yielding optimal control sequence $u_i^*(k|t), k = 0,1,\cdots,N_p - 1$.
2) Compute optimal state in the predictive horizon using optimal control $u_i^*(k|t)$

$$x_i^*(k+1|t) = \phi_i\bigl(x_i^*(k|t)\bigr) + \psi_i \cdot u_i^*(k|t),$$
$$k = 0,1,\cdots,N_p - 1. \quad (16)$$
$$x_i^*(0|t) = x_i(t).$$

3) Compute the assumed control input (*i.e.*, $u_i^a(k|t+1)$) for next step by disposing first term and adding one additional term, *i.e.*,

$$u_i^a(k|t+1) = \begin{cases} u_i^*(k+1|t), & k = 0,1,\cdots,N_p - 2 \\ h_i\bigl(v_i^*(N_p|t)\bigr), & k = N_p - 1 \end{cases}. \quad (17)$$

The corresponding assumed output is also computed as

$$x_i^a(k+1|t+1) = \phi_i\bigl(x_i^a(k|t+1)\bigr)$$
$$+ \psi_i u_i^a(k|t+1),$$
$$x_i^a(0|t+1) = x_i^*(1|t), \quad (18)$$
$$y_i^a(k|t+1) = \gamma x_i^a(k|t+1),$$
$$k = 0,1,\cdots,N_p - 1.$$

4) Transmit $y_i^a(k|t+1)$ to the nodes in set $\mathbb{O}_i$, receive $y_{-i}^a(k|t+1)$ from the nodes in set $\mathbb{N}_i$, and then compute $y_{\text{des},i}(k|t+1)$ using the leader's information if $\mathbb{P}_i \neq \emptyset$.
5) Implement the control effort using the first element of optimal control sequence, *i.e.*, $u_i(t) = u_i^*(0|t)$.
6) Increment $t$ and go to step (1).

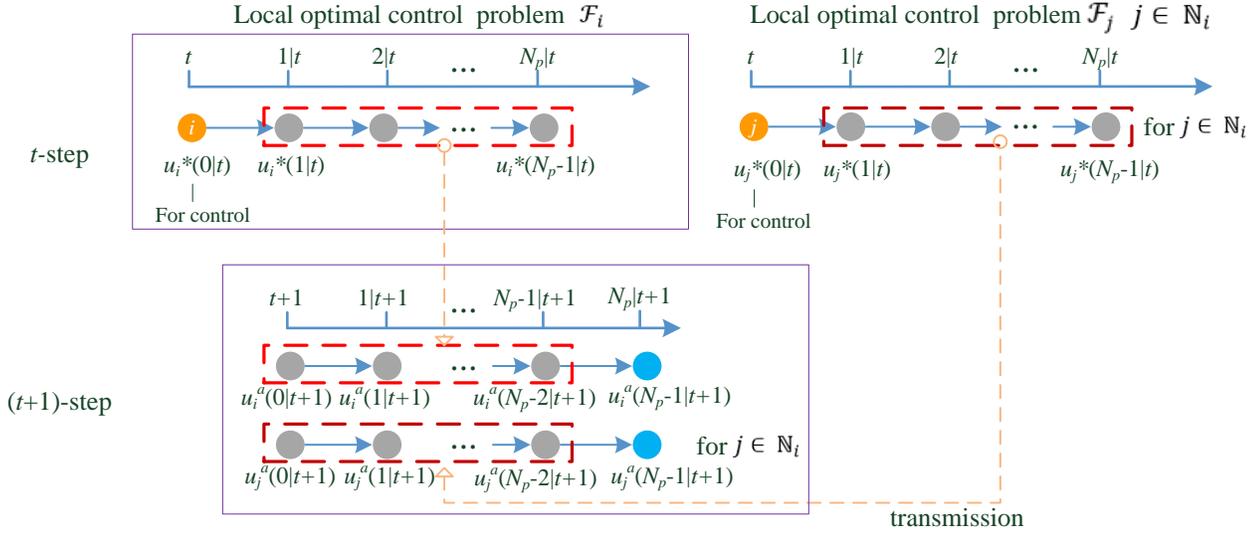

Fig. 3. The basic procedure to construct assumed inputs.

**Remark 4**. One key part of DMPC is how to construct the assumed input and output in each node. Here, the assumed variable is a shifted optimal result of last-step problem $\mathcal{F}_i$, synthesized by disposing the first value and adding a last value. The last added value ensures that the vehicle moves at a constant speed. A similar technique can be found in [14] and [27]. Fig. 3 gives schematic procedure to construct assumed inputs. Note that in this DMPC framework, all followers are assumed to be synchronized in the step of control execution, *i.e.*, updating the system state simultaneously within a common global clock. However, neither computation nor communication is assumed to happen instantaneously.

**Remark 5**. Another key feature of this DMPC algorithm is that each node only needs to solve a local optimization problem of small size relying on the information of its neighbors in set $\mathbb{N}_i$, and pass the results to the nodes in set $\mathbb{O}_i$ at each time step. Additionally, the computational complexity of $\mathcal{F}_i$ is independent with the platoon size $N$, which implies the proposed DMPC approach is scalable provided a single MPC in each node can be solved efficiently. In this aspect, several efficient computing techniques, such as utilizing particular structure [40], using explicit MPC via a lookup table [41], and reducing the dimension via the parameterization method and "move blocking" method [42], might be employed to solve each single MPC problem for real-time implementations, which would be extremely interesting for further research.

## IV. STABILITY ANALYSIS OF THE DMPC ALGORITHM

This section presents the stability analysis of the proposed DMPC algorithm. The main strategy is to construct a proper Lyapunov candidate for the platoon and prove its decreasing property. A sufficient condition for asymptotic stability is derived by using the sum of local cost functions as a Lyapunov function. The condition shows that stability can be achieved through explicit sufficient conditions on the weights in the cost functions.

### A. Terminal constraint analysis

For the completeness of proof, we first present the assumption made on allowable topologies.

*Assumption 1*: The graph $\mathbb{G}$ contains a spanning tree rooting at the leader, and the communications are unidirectional from preceding vehicles to downstream ones.

The topologies satisfying abovementioned **Assumption 1** are called unidirectional topology for short. Fig. 1 shows some typical examples. The following lemmas are useful for stability analysis.

*Lemma 1* [35]. Suppose that $\lambda_1, \ldots, \lambda_n$ be the eigenvalues of $A \in \mathbb{R}^{n \times n}$ and $\mu_1, \ldots, \mu_m$ be those of $B \in \mathbb{R}^{m \times m}$. Then, the eigenvalues of $A \otimes B$ are
$$\lambda_i \mu_j, i = 1, \cdots, n, j = 1, \cdots, m.$$

*Lemma 2* [35]. Let a matrix $Q = [q_{ij}] \in \mathbb{R}^{n \times n}$. Then, all the eigenvalues of $Q$ are located in the union of the $n$ disks
$$\cup_{i=1}^n \{\lambda \in \mathbb{C} | |\lambda - q_{ii}| \leq \sum_{j=1, j \neq i}^n |q_{ij}|\}.$$

This is the well-known Geršgorin Disk Criterion.

*Lemma 3* [21][36]. Matrix $\mathcal{L} + \mathcal{P}$ is nonsingular if $\mathbb{G}$ contains a spanning tree rooting at the leader.

*Lemma 4*. If $\mathbb{G}$ contains a spanning tree rooting at the leader, then $\mathcal{D} + \mathcal{P}$ is invertible and all the eigenvalues of $(\mathcal{D} + \mathcal{P})^{-1} \mathcal{A}$ are located within a unit circle, *i.e.*,
$$\{\lambda \in \mathbb{C} | |\lambda| < 1\} \quad (19)$$

**Proof.** Since $\mathbb{G}$ contains a spanning tree, for any $i \in \mathcal{N}$, we have that either the in-degree of node $i$ is larger than zero, *i.e.*, $\deg_i > 0$, or node $i$ is pinned to the leader, *i.e.*, $p_i = 1$, or both of them are true. Either way, we know $\deg_i + p_i \geq 1$, and considering the fact that $\mathcal{D} + \mathcal{P}$ is a diagonal matrix, we have
$$\mathcal{D} + \mathcal{P} > 0. \quad (20)$$

Thus, $\mathcal{D} + \mathcal{P}$ is invertible.

Let $\sigma_i, i \in \mathcal{N}$ be the eigenvalues of $(\mathcal{D} + \mathcal{P})^{-1} \mathcal{A}$. Considering the definition (7), the diagonal elements of

$(\mathcal{D} + \mathcal{P})^{-1}\mathcal{A}$ are all equal to zero. Then, according to **Lemma 2**, $\sigma_i, i \in \mathcal{N}$ are located in the union of $N$ disks:

$$\bigcup_{i=1}^{N} \left\{ \lambda \in \mathbb{C} | |\lambda - 0| \leq \sum_{j=1, j \neq i}^{N} \left| \frac{a_{ij}}{\deg_i + p_i} \right| \right\}. \quad (21)$$

Further, we have

$$\sum_{j=1, j \neq i}^{N} \left| \frac{a_{ij}}{\deg_i + p_i} \right| = \left| \frac{\deg_i}{\deg_i + p_i} \right| \leq 1. \quad (22)$$

Combining (21) and (22), $\sigma_i, i \in \mathcal{N}$ are bounded by a unit circle

$$\{\lambda \in \mathbb{C} | |\lambda| \leq 1\}. \quad (23)$$

Next we will prove that $\sigma_i$ cannot be located on the boundary of the unit circle by contradiction. Suppose that some eigenvalues are located on the boundary, *i.e.*,

$$\rho((\mathcal{D} + \mathcal{P})^{-1}\mathcal{A}) = 1. \quad (24)$$

Since $(\mathcal{D} + \mathcal{P})^{-1}\mathcal{A}$ is non-negative, one of its eigenvalues is equal to one according to (24) [35]. Let the corresponding eigenvector be $x$, then the following equality holds.

$$(\mathcal{D} + \mathcal{P})^{-1}\mathcal{A} \cdot x = x. \quad (25)$$

Considering the fact $\mathcal{A} = \mathcal{D} - \mathcal{L}$, we have

$$(\mathcal{L} + \mathcal{P}) \cdot x = 0. \quad (26)$$

Then, $\mathcal{L} + \mathcal{P}$ is singular from (26), which is in contradiction with **Lemma 3**. Therefore, $\sigma_i$ cannot be located on the boundary, which means

$$\{\lambda \in \mathbb{C} | |\lambda| < 1\}. \quad (27)$$

∎

Here, we have the following theorem.

***Theorem 1***. If $\mathbb{G}$ contains a spanning tree rooting at the leader, the terminal state in the predictive horizon of problem $\mathcal{F}_i$ asymptotically converges to the desired state, *i.e.*,

$$\lim_{t \to \infty} |y_i^p(N_p|t) - y_{\text{des},i}(N_p|t)| = 0. \quad (28)$$

where $y_{\text{des},i}(N_p|t) = [s_0(N_p|t) - i \cdot d_0, v_0]^T$.

**Proof**. Constrained by (13e), each node moves at constant speed at the end of predictive horizon. Considering assumed control input (17), we have

$$y_i^a(N_p|t+1) = y_i^p(N_p|t) + B \cdot y_i^p(N_p|t) \cdot \Delta t,$$
$$B = \begin{bmatrix} 0 & 1 \\ 0 & 0 \end{bmatrix}. \quad (29)$$

Submitting (29) into (13d) yields

$$y_i^p(N_p|t+1) = \frac{1}{|\mathbb{I}_i|} \sum_{j \in \mathbb{I}_i} (y_j^p(N_p|t) + B y_j^p(N_p|t) \\
\cdot \Delta t + \tilde{d}_{i,j}). \quad (30)$$

Define the tracking error vector as

$$\hat{y}_i^p(N_p|t) = y_i^p(N_p|t) - y_{\text{des},i}(N_p|t), \quad (31)$$

and we have (32) by combining (30) and (31)

$$\hat{y}_i^p(N_p|t+1) = \frac{1}{|\mathbb{I}_i|} \sum_{j \in \mathbb{I}_i} (I_2 + B\Delta t) \hat{y}_j^p(N_p|t). \quad (32)$$

Define the collected terminal state vector as $Y^p(N_p|t) = [\hat{y}_j^p(N_p|t), \cdots, \hat{y}_j^p(N_p|t)]^T \in \mathbb{R}^{2N \times 1}$, then (32) can be further written into a compact form

$$Y^p(N_p|t+1) = [(\mathcal{D} + \mathcal{P})^{-1} \cdot \mathcal{A}] \otimes (I_2 + B\Delta t) \\
\cdot Y^p(N_p|t). \quad (33)$$

It is easy to verify that the eigenvalues of $I_2 + B\Delta t$ are all equal to one. Besides, according to **Lemma 4**, all the eigenvalues of $(\mathcal{D} + \mathcal{P})^{-1}\mathcal{A}$ are located within a unit circle. Thus, by **Lemma 1**, the eigenvalues of $[(\mathcal{D} + \mathcal{P})^{-1} \cdot \mathcal{A}] \otimes (I_2 + B\Delta t)$ are all located within a unit circle as well, *i.e.*,

$$\{\lambda \in \mathbb{C} | |\lambda| < 1\}. \quad (34)$$

Then, based on (33), we know $Y^p(N_p|t)$ asymptotically converges to zero, which means

$$\lim_{t \to \infty} |y_i^p(N_p|t) - y_{\text{des},i}(N_p|t)| = 0. \quad (35)$$

∎

***Theorem 2***. If $\mathbb{G}$ satisfies **Assumption 1**, the terminal state in the predictive horizon of problem $\mathcal{F}_i$ converges to the desired state in at most $N$ steps, *i.e.*,

$$y_i^p(N_p|t) = y_{\text{des},i}(N_p|t), \quad t \geq N. \quad (36)$$

**Proof**. If $\mathbb{G}$ is unidirectional, then $\mathcal{A}$ is a lower triangular matrix with diagonal entries be zero. Based on (20), $\mathcal{D} + \mathcal{P} > 0$. Therefore, the eigenvalues of $(\mathcal{D} + \mathcal{P})^{-1}\mathcal{A}$ are all zero, and $(\mathcal{D} + \mathcal{P})^{-1}\mathcal{A}$ is nilpotent with degree at most $N$.

By **Lemma 1**, we further have that the eigenvalues of $[(\mathcal{D} + \mathcal{P})^{-1} \cdot \mathcal{A}] \otimes (I_2 + B\Delta t)$ are all zero as well. Hence, $Y^p(N_p|t)$ can converge to zero in at most $N$ steps, which means that $y_i^p(N_p|t)$ in $\mathcal{F}_i$ converges to the desired state in at most $N$ steps.

∎

**Remark 6.** Even though not every follower directly communicates with the leader, the terminal state of each node can still converge to its desired set point within finite time under **Assumption 1**, which means this DMPC scheme does not require all nodes *a priori* know the desired set point. Note that in the proposed DMPC algorithm, the number of time steps required for the consensus of the terminal states is upper bounded by the platoon size (*i.e.*, $N$). This implies an intuitive fact that the speed of sharing leader's information is directly affected by the size of a platoon for unidirectional topologies.

**Remark 7.** For homogenous platoons with linear dynamics and linear controllers, it is well demonstrated that stability requires at least a spanning tree rooting at the leader [21][36]. Even using an MPC technique, a spanning tree is also a prerequisite to achieve a stable platoon. Intuitively, this requirement means that every follower can obtain the leader information directly or indirectly.

**Remark 8.** The length of predictive horizon $N_p$ has no explicit relationship with platoon size $N$ in terms of asymptotic stability. It should be note that the analysis of terminal constraint relies on the assumption that each local optimization problem $\mathcal{F}_i$ is feasible for the first $N$ steps. This is called initial feasible assumption, which is widely used in previous studies on DMPC, *e.g.*, [14][27]-[29][39]. After the consensus of the terminal states, the property of recursive feasibility holds (see **Lemma 5**). Consequently, $N_p$ should be large enough to get a feasible solution for problem $\mathcal{F}_i$ (note that initial errors will also affect the feasibility in addition to the model and constraints). However, a large $N_p$ will lead to a great computing burden in terms of computation time and memory requirement. The optimal choice of time horizon $N_p$ should be a balance of performance and computational effort [37], which is beyond the scope of this paper.

### B. Analysis of local cost function

The optimal cost function of node $i$ at time $t$ is denoted as

$$J_i^*(t) = J_i^*\big(y_i^*(:|t), u_i^*(:|t), y_i^a(:|t), y_{-i}^a(:|t)\big). \tag{37}$$

The following is a standard result in MPC formulation.

**Lemma 5** [14]. If we replace (13d) with $y_i^p(N_p|t) = y_{\text{des},i}(N_p|t)$, then problem $\mathcal{F}_i$ has

$$\big(y_i^p(:|t), u_i^p(:|t)\big) = \big(y_i^a(:|t), u_i^a(:|t)\big) \tag{38}$$

as a feasible solution for any time $t > 0$.

Note that **Lemma 5** is the property of recursive feasibility. The assumed control $u_i^a(:|t)$ defined in (17) is the same feasible control used in [14] and [27]. The remaining part of this section is to analyze the decreasing properties of local cost function. Here, we have the following theorem.

**Theorem 3**. If $\mathbb{G}$ satisfies **Assumption 1**, each local cost function satisfies

$$\begin{aligned}&J_i^*(t+1) - J_i^*(t)\\ &\leq -l_i\big(y_i^*(0|t), u_i^*(0|t), y_i^a(0|t), y_{-i}^a(0|t)\big) + \varepsilon_i,\\ &\qquad t > N\end{aligned} \tag{39}$$

where

$$\varepsilon_i = \sum_{k=1}^{N_p-1} \left\{ \sum_{j\in\mathbb{N}_i} \big\|y_j^*(k|t) - y_j^a(k|t)\big\|_{G_i} \right.$$
$$\left. - \big\|y_i^*(k|t) - y_i^a(k|t)\big\|_{F_i} \right\}.$$

**Proof**. If $\mathbb{G}$ satisfies **Assumption 1**, **Theorem 2** gives that $y_i^p(N_p|t) - y_{\text{des},i}(N_p|t) = 0$, $t \geq N$. Then, at time $t+1, t \geq N$, a feasible (but suboptimal) control for $\mathcal{F}_i$ is $u_i^p(:|t+1) = u_i^a(:|t+1)$. Therefore, we can bound the optimal cost as

$$\begin{aligned}J_i^*(t+1) &\leq J_i\big(y_i^a(:|t+1), u_i^a(:|t+1), y_i^a(:|t+1), y_{-i}^a(:|t+1)\big)\\ &= \sum_{k=0}^{N_p-1} l_i\big(y_i^a(k|t+1), u_i^a(k|t+1), y_i^a(k|t+1), y_{-i}^a(k|t+1)\big)\end{aligned} \tag{40}$$

$$= \sum_{k=0}^{N_p-2} l_i\big(y_i^*(k+1|t), u_i^*(k+1|t), y_i^*(k+1|t), y_{-i}^*(k+1|t)\big).$$

The equality holds because of how $u_i^a(k|t+1)$ and $y_i^a(k|t+1)$ are defined by (17) and (18).

Further, by changing the index of summation, (40) becomes

$$\begin{aligned}&J_i^*(t+1)\\ &\leq \sum_{k=1}^{N_p-1} l_i\big(y_i^*(k|t), u_i^*(k|t), y_i^*(k|t), y_{-i}^*(k|t)\big).\end{aligned} \tag{41}$$

Subtracting $J_i^*(t)$ from (41) yields

$$\begin{aligned}&J_i^*(t+1) - J_i^*(t)\\ &\leq \sum_{k=1}^{N_p-1} l_i\big(y_i^*(k|t), u_i^*(k|t), y_i^*(k|t), y_{-i}^*(k|t)\big)\\ &\quad - \sum_{k=0}^{N_p-1} l_i\big(y_i^*(k|t), u_i^*(k|t), y_i^a(k|t), y_{-i}^a(k|t)\big)\\ &= -l_i\big(y_i^*(0|t), u_i^*(0|t), y_i^a(0|t), y_{-i}^a(0|t)\big) + \sum_{k=1}^{N_p-1}\Delta_k,\end{aligned} \tag{42}$$

where

$$\begin{aligned}\Delta_k &= l_i\big(y_i^*(k|t), u_i^*(k|t), y_i^*(k|t), y_{-i}^*(k|t)\big)\\ &\quad - l_i\big(y_i^*(k|t), u_i^*(k|t), y_i^a(k|t), y_{-i}^a(k|t)\big)\\ &= \big\|y_i^*(k|t) - y_{\text{des},i}(k|t)\big\|_{Q_i}\\ &\quad + \big\|u_i^*(k|t) - h_i\big(v_i^*(k|t)\big)\big\|_{R_i}\\ &\quad + \big\|y_i^*(k|t) - y_i^*(k|t)\big\|_{F_i}\\ &\quad + \sum_{j\in\mathbb{N}_i}\big\|y_i^*(k|t) - y_j^*(k|t) - \tilde{d}_{i,j}\big\|_{G_i}\\ &\quad - \Big\{\big\|y_i^*(k|t) - y_{\text{des},i}(k|t)\big\|_{Q_i} + \big\|u_i^*(k|t) - h_i\big(v_i^*(k|t)\big)\big\|_{R_i}\\ &\quad + \big\|y_i^*(k|t) - y_i^a(k|t)\big\|_{F_i}\\ &\quad + \sum_{j\in\mathbb{N}_i}\big\|y_i^*(k|t) - y_j^a(k|t) - \tilde{d}_{i,j}\big\|_{G_i}\Big\}.\end{aligned} \tag{43}$$

With the triangle inequality for vector norms, (43) becomes

$$\begin{aligned}\Delta_k &\leq \sum_{j\in\mathbb{N}_i}\big\|y_j^*(k|t) - y_j^a(k|t)\big\|_{G_i}\\ &\quad - \big\|y_i^*(k|t) - y_i^a(k|t)\big\|_{F_i}.\end{aligned} \tag{44}$$

Combining (42) and (44) yields (39). ∎

**Remark 9.** Note that (39) gives an upper bound on the decline of local cost function. If we have

$$\varepsilon_i \leq l_i\big(y_i^*(0|t), u_i^*(0|t), y_i^a(0|t), y_{-i}^a(0|t)\big), \tag{45}$$

then local cost function decreases monotonically, which means it is a proper Lyapunov function and also leads to asymptotic stability of DMPC. The difficulty of using (45) to design DMPC is obvious, *i.e.*, there is no intuitive way to adjust control parameters. One alternative is to use the sum of local cost functions as a Lyapunov function as suggested by [28].

### C. Sum of local cost functions

Define the sum of all local cost functions as the Lyapunov candidate:

$$J_\Sigma^*(t) = \sum_{i=1}^{N} J_i^*\big(y_i^*(:|t), u_i^*(:|t), y_i^a(:|t), y_{-i}^a(:|t)\big). \quad (46)$$

Then, we have the following theorem.

***Theorem 4***. If $\mathbb{G}$ satisfies **Assumption 1**, $J_\Sigma^*(t)$ satisfies

$$\begin{aligned} J_\Sigma^*(t+1) - J_\Sigma^*(t) \\ \leq -\sum_{i=1}^{N} l_i\big(y_i^*(0|t), u_i^*(0|t), y_i^a(0|t), y_{-i}^a(0|t)\big) \\ + \sum_{k=1}^{N_p-1} \varepsilon_\Sigma(k), \ t > N \end{aligned} \quad (47)$$

where

$$\varepsilon_\Sigma(k) = \sum_{i=1}^{N} \left[ \sum_{j \in \mathbb{O}_i} \|y_i^*(k|t) - y_i^a(k|t)\|_{G_j} \\ - \|y_i^*(k|t) - y_i^a(k|t)\|_{F_i} \right].$$

**Proof**. According to **Theorem 3**, we have

$$\begin{aligned} J_\Sigma^*(t+1) - J_\Sigma^*(t) \\ \leq \sum_{i=1}^{N} \big\{-l_i\big(y_i^*(0|t), u_i^*(0|t), y_i^a(0|t), y_{-i}^a(0|t)\big) + \varepsilon_i\big\} \\ = -\sum_{i=1}^{N} l_i\big(y_i^*(0|t), u_i^*(0|t), y_i^a(0|t), y_{-i}^a(0|t)\big) \\ + \sum_{i=1}^{N} \varepsilon_i. \end{aligned} \quad (48)$$

Further, we know

$$\begin{aligned} \sum_{i=1}^{N} \varepsilon_i = \sum_{k=1}^{N_p-1} \Bigg\{ \sum_{i=1}^{N} \left[ \sum_{j \in \mathbb{N}_i} \|y_j^*(k|t) - y_j^a(k|t)\|_{G_i} \\ - \|y_i^*(k|t) - y_i^a(k|t)\|_{F_i} \right] \Bigg\} \\ = \sum_{k=1}^{N_p-1} \Bigg\{ \sum_{i=1}^{N} \left[ \sum_{j \in \mathbb{O}_i} \|y_i^*(k|t) - y_i^a(k|t)\|_{G_j} \\ - \|y_i^*(k|t) - y_i^a(k|t)\|_{F_i} \right] \Bigg\} \end{aligned} \quad (49)$$

$$= \sum_{k=1}^{N_p-1} \varepsilon_\Sigma(k).$$

Combining (48) and (49) yields (47). ∎

**Remark 10.** The key in the proof of **Theorem 4** is to change $\mathbb{N}_i$ to $\mathbb{O}_i$ by considering all followers in the platoon. Note that (47) is also an upper bound on the decline of the sum of local cost function. Moreover, it is relatively easy for designers to find a sufficient condition to guarantee $J_\Sigma^*(t+1) - J_\Sigma^*(t) < 0$.

### D. Sufficient condition of DMPC stability

The explicit sufficient stability condition is now stated as follows.

***Theorem 5***. If $\mathbb{G}$ satisfies **Assumption 1**, a platoon under DMPC (13) is asymptotically stable if satisfying

$$F_i \geq \sum_{j \in \mathbb{O}_i} G_j, \ i \in \mathcal{N}. \quad (50)$$

**Proof**:
If (50) holds, we have

$$z^T \left( \sum_{j \in \mathbb{O}_i} G_j - F_i \right) z \leq 0, \forall z \in \mathbb{R}^2 \quad (51)$$

Let $z = y_i^*(k|t) - y_i^a(k|t)$, then

$$\sum_{j \in \mathbb{O}_i} \|y_i^*(k|t) - y_i^a(k|t)\|_{G_j} \\ - \|y_i^*(k|t) - y_i^a(k|t)\|_{F_i} \leq 0. \quad (52)$$

Combining **Theorem 4**, we have

$$\begin{aligned} J_\Sigma^*(t+1) - J_\Sigma^*(t) \\ \leq -\sum_{i=1}^{N} l_i\big(y_i^*(0|t), u_i^*(0|t), y_i^a(0|t), y_{-i}^a(0|t)\big). \end{aligned} \quad (53)$$

The upper bound in (53) shows that $J_\Sigma^*(t)$ is strictly monotonically decreasing. Thus, the asymptotic stability of DMPC is guaranteed. ∎

**Remark 11. Theorem 5** shows that for heterogeneous platoons under unidirectional topologies, it only needs to adjust the weights on the errors between predicted trajectories and assumed ones to guarantee asymptotic stability. Note that the condition (50) in **Theorem 5** is distributed with respect to the vehicles in the platoon. The followers in a platoon do not need the centralized information to choose their own penalty weights.

**Remark 12.** Notation $\mathbb{O}_i = \{j|a_{ji} = 1\}$ in (50) is defined as the nodes that can use the information of node $i$. This provides an interesting phenomenon, *i.e.*, to ensure stability implies that all nodes in $\mathbb{O}_i$ should not rely heavily on the information of node $i$ unless node $i$ shows good-enough consistence with its own assumed trajectory.

**Remark 13.** Many previous studies on platoon control using DMPC techniques are only suitable for some special topologies, for example [14] and [34]. **Theorem 5** extends the topological selection to be any arbitrary unidirectional topology (defined in **Assumption 1**), which can include many other types of topologies (see Fig. 1 for examples).

## V. SIMULATION RESULTS

In this section, numerical simulations are conducted to illustrate the main results of this paper. We consider a heterogeneous platoon with eight vehicles (*i.e.*, 1 leader and 7 followers) interconnected by the four types of communication topologies shown in Fig. 1.

The acceleration of the leader can be viewed as disturbances in a platoon [15][22]. The initial state of the leader is set as $s_0(t) = 0, v_0 = 20$ m/s and the desired trajectory is given by

$$v_0 = \begin{cases} 20 \text{ m/s} & t \leq 1 \text{ s} \\ 20 + 2t \text{ m/s} & 1s < t \leq 2 \text{ s.} \\ 22 \text{ m/s} & t > 2s \end{cases}$$

The parameters of following vehicles are randomly selected according to the passenger vehicles [18], which are listed in TABLE I. In the simulation, the box constraints on tracking/braking torque are reflected by the maximum acceleration and deceleration, *i.e.*, $a_{max,i} = 6 \text{ m/s}^2$, $a_{min,i} = -6 \text{ m/s}^2$. The discrete time interval is chosen as $\Delta t = 0.1$ s, and the predictive horizon in the $\mathcal{F}_i$ is all set as $N_p = 20$. TABLE II lists the corresponding weights in the $\mathcal{F}_i$, which can be easily verified to satisfy the conditions in **Theorem 5**.

In the simulations, the desired spacing is set as $d_{i-1,i} = 20$ m. The initial state of the platoon is set as the desired state, *i.e.,* the initial spacing errors and velocity errors are all equal to 0. Fig. 4 demonstrates the spacing errors of the platoon under different topologies. It is easy to find that the platoon using DMPC is stable for all topologies listed in Fig. 1, which conforms to the results in **Theorem 5**. Additionally, for this simulation scenario, the spacing errors are less than 1 m for platoons with all of the four communication topologies. This result also shows that there are no collisions during the transient process.

TABLE I. PARAMETERS OF THE FOLLOWING VEHICLES IN THE PLATOON

| Vehicle Index | $m_{\text{veh},i}$ (kg) | $\tau_i$ (s) | $C_{A,i}$ (N·s²·m⁻²) | $R_i$ (m) |
|---|---|---|---|---|
| 1 | 1035.7 | 0.51 | 0.99 | 0.30 |
| 2 | 1849.1 | 0.75 | 1.15 | 0.38 |
| 3 | 1934.0 | 0.78 | 1.17 | 0.39 |
| 4 | 1678.7 | 0.70 | 1.12 | 0.37 |
| 5 | 1757.7 | 0.73 | 1.13 | 0.38 |
| 6 | 1743.1 | 0.72 | 1.13 | 0.37 |
| 7 | 1392.2 | 0.62 | 1.06 | 0.34 |

TABLE II. WEIGHTS IN THE COST FUNCTONS

| Weights | PF | PLF | TPF | TPLF |
|---|---|---|---|---|
| $F_i$ | $F_i = 10I_2$, $i \in \mathcal{N}$ | $F_i = 10I_2$, $i \in \mathcal{N}$ | $F_i = 10I_2$, $i \in \mathcal{N}$ | $F_i = 10I_2$, $i \in \mathcal{N}$ |
| $G_i$ | $G_1 = 0$, $G_i = 5I_2$, $i \in \mathcal{N}\backslash\{1\}$ | $G_1 = 0$, $G_i = 5I_2$, $i \in \mathcal{N}\backslash\{1\}$ | $G_1 = 0$, $G_i = 5I_2$, $i \in \mathcal{N}\backslash\{1\}$ | $G_1 = 0$, $G_i = 5I_2$, $i \in \mathcal{N}\backslash\{1\}$ |
| $Q_i$ | $Q_1 = 10I_2$, $Q_i = 0$, $i \in \mathcal{N}\backslash\{1\}$ | $Q_i = 10I_2$, $i \in \mathcal{N}$ | $Q_1 = 10I_2$, $Q_2 = 10I_2$, $Q_i = 0$, $i \in \mathcal{N}\backslash\{1,2\}$ | $Q_i = 10I_2$, $i \in \mathcal{N}$ |
| $R_i$ | $R_i = I_2$, $i \in \mathcal{N}$ | $R_i = I_2$, $i \in \mathcal{N}$ | $R_i = I_2$, $i \in \mathcal{N}$ | $R_i = I_2$, $i \in \mathcal{N}$ |

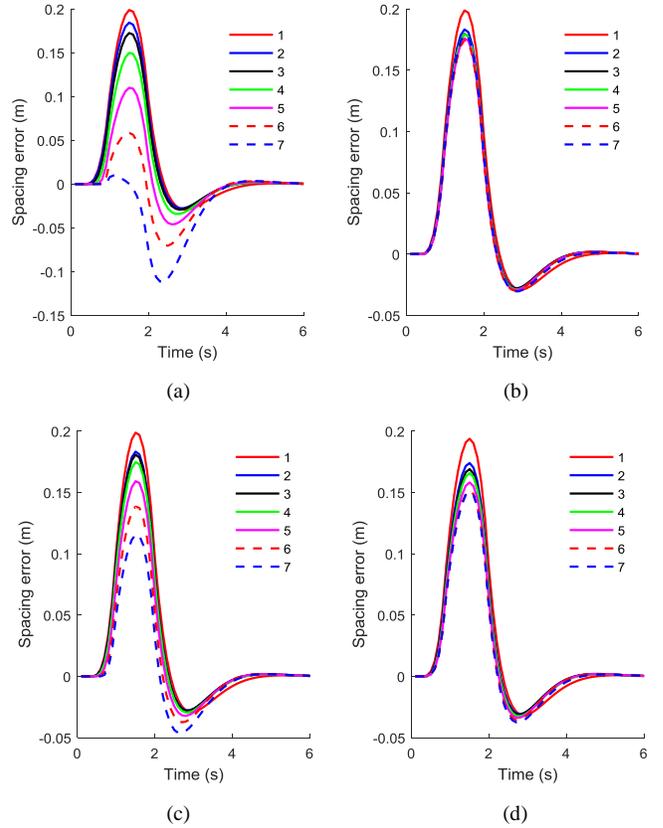

Fig. 4. Spacing errors for the platoon under different topologies. (a) PF; (b) PLF; (c) TPF; (d) TPLF;

## VI. CONCLUSIONS

This paper proposes a novel DMPC algorithm for vehicle platoons with nonlinear dynamics and unidirectional topologies, and derives a sufficient condition to guarantee asymptotic stability. This approach does not require all vehicle *a priori* know the desired set point, which offers considerable benefit from the viewpoint of real implementations.

Under the proposed DMPC framework, the platoon is dynamically decoupled, but constrained by the spatial formation. Each vehicle has nonlinear dynamics with input constraints, but does not necessarily know its desired set point. Each vehicle solves a local optimal control problem to obtain its own control input, and then sends its assumed output trajectory to its neighbors. A neighboring average based terminal constraint is introduced in the formulation of local optimal problems, which guarantees that all terminal states in the predictive horizon can converge to the desired state in finite time when the topology is unidirectional and contains a spanning tree. By using the sum of local cost functions as the Lyapunov function, it is further proved that asymptotic stability can be achieved through an explicit sufficient condition on the weights of the cost functions.

One topic for future research is to improve DMPC algorithm by deriving the stability condition under more general topologies. Besides, other important issues include how to address the disturbances and uncertainty in the dynamics, and how to handle the packets drops and delays in the communication between vehicles.